\documentclass[a4paper,11pt]{amsart}
\usepackage{graphicx}
\usepackage{mathptmx}
\usepackage{amsmath}
\usepackage{amssymb}
\usepackage{enumitem}
\usepackage{xcolor}

\newtheorem{theorem}{Theorem}[section]
\newtheorem{lemma}[theorem]{Lemma}

\newcommand{\sech}{\text{sech}}
\newcommand{\csch}{\text{csch}}
\newcommand{\intz}{\int_{\mathbb{Z}_p}}
\newcommand{\dm}{d \mu_1}
\newcommand{\dmm}{d \mu_{-1}}

\begin{document}

\title{Some identities of special numbers and polynomials arising from $p$-adic integrals on $\mathbb{Z}_p$}

\author{Dae San Kim}
\address{Department of Mathematics, Sogang University, Seoul 04107, Republic of Korea}
\email{dskim@sogang.ac.kr}

\author{Han Young Kim}
\address{Department of Mathematics, Kwangwoon University, Seoul 01897, Republic of Korea}
\email{gksdud213@gmail.com}

\author{Sung-Soo Pyo}
\address{Department of Mathematics Education, Silla University, Busan 46958, Republic of Korea}
\email{ssoopyo@gmail.com}

\author{Taekyun Kim}
\address{Department of Mathematics, Kwangwoon University, Seoul 01897, Republic of Korea}
\email{tkkim@kw.ac.kr}

\maketitle

\begin{abstract}
In recent years, studying degenerate versions of various special polynomials and numbers have attracted many mathematicians. Here we introduce degenerate type 2 Bernoulli polynomials, fully degenerate type 2 Bernoulli polynomials and degenerate type 2 Euler polynomials, and their corresponding numbers, as degenerate and type 2 versions of Bernoulli and Euler numbers. Regarding to those polynomials and numbers, we derive some identities, distribution relations, Witt type formulas and analogues for the Bernoulli's interpretation of powers of the first $m$ positive integers in terms of Bernoulli polynomials. The present study was done by using the bosonic and fermionic $p$-adic integrals on $\mathbb{Z}_p$.
\end{abstract}

\markboth{\centerline{\scriptsize Some identities of special numbers and polynomials arising from $p$-adic integrals on $\mathbb{Z}_p$ }}
{\centerline{\scriptsize D. S. Kim, H. Y. Kim, S.-S. Pyo, T. Kim}}

\section{Introduction}\label{sec1}

\indent Studies on degenerate versions of some special polynomials and numbers began with the papers by Carlitz in \cite{ref3, ref4}. In recent years, studying degenerate versions of various special polynomials and numbers have regained interests of many mathematicians. The researches have been carried out by several different methods like generating functions, combinatorial approaches, umbral calculus, $p$-adic analysis and differential equations. This idea of studying degenerate versions of some special polynomials and numbers turned out to be very fruitful so as to introduce degenerate Laplace transforms and degenerate gamma functions (see~\cite{ref12}).\\
\indent In this paper, we introduce degenerate type 2 Bernoulli polynomials, fully degenerate type 2 Bernoulli polynomials and degenerate type 2 Euler polynomials, and their corresponding numbers, as degenerate and type 2 versions of Bernoulli and Euler numbers. We investigate those polynomials and numbers by means of  bosonic and fermionic $p$-adic integrals and derive some identities, distribution relations, Witt type formulas and analogues for the Bernoulli's interpretation of powers of the first $m$ positive integers in terms of Bernoulli polynomials. In more detail, our main results are as follows.\\
\indent As to the analogues for the Bernoulli's interpretation of power sums, in Theorem 2.6 we express powers of the first $m$ odd integers in terms of type 2 Bernoulli polynomials $b_n(x)$, in Theorem 2.11  alternating sum of powers of the first $m$ odd integers in terms of type 2 Euler polynomials $E_n(x)$, in Theorem 2.9 sum of the values of the generalized falling factorials at the first $m$ odd positive integers in terms of degenerate Carlitz type 2 Bernoulli polynomials $b_{n,\lambda}(x)$, and in Theorem 2.17 alternating sum of the values of the generalized falling factorials at the first $m$ odd positive integers in terms of degenerate type 2 Euler polynomials $E_{n,\lambda}(x)$. Witt type formulas are obtained for $b_n(x), B_{n,\lambda}(x), E_n(x)$, and $E_{n,\lambda}(x)$, respectively in Lemma 2.1, Theorem 2.7, Lemma 2.10 and Theorem 2.16. 
Distribution relations are derived for  $b_n(x)$, and $E_n(x)$, respectively in Theorem 2.3 and Theorem 2.13. \\
\indent In the rest of this section, we will introduce type 2 Bernoulli and Euler numbers, recall the bosonic and fermionic $p$-adic integrals and mention the degenerate exponential function.

Let $p$ be a fixed odd prime number. Throughout this paper, $\mathbb{Z}_p$, $\mathbb{Q}_p$ and $\mathbb{C}_p$ will denote the ring of $p$-adic integers, the field of $p$-adic rational numbers and the completion of an algebraic closure of $\mathbb{Q}_p$, respectively.
The $p$-adic norm $|\cdot|_p$ is normalized by $|p|_p=\frac{1}{p}$.\\
\indent It is well known that the ordinary Bernoulli polynomials are defined by 
\begin{equation}\label{eq1}
\frac{t}{e^t-1}e^{xt}=\sum_{n=0}^\infty B_n(x)\frac{t^n}{n!},\quad (\text{see \cite{ref2, ref5, ref14, ref15, ref17}}).
\end{equation}
When $x=0$, $B_n=B_n(0)$ are called the Bernoulli numbers. 

Also, the type 2 Bernoulli polynomials are given by 
\begin{equation}\label{eq2}
\frac{t}{2}\csch\frac{t}{2}
e^{xt}=\frac{t}{e^{\frac{t}{2}}-e^{-\frac{t}{2}}}e^{xt}
=\sum_{n=0}^\infty b_n(x)\frac{t^n}{n!}.
\end{equation}
For $x=0$, $b_n=b_n(0)$ are called the type 2 Bernoulli numbers so that they are given by
\begin{equation}\label{eq3}
\frac{t}{2}\csch
\frac{t}{2}=\frac{t}{e^{\frac{t}{2}}-e^{-\frac{t}{2}}}
=\sum_{n=0}^\infty b_n\frac{t^n}{n!}.
\end{equation}
In fact, the type 2 Bernoulli polynomials and numbers are slightly differently defined in \cite{ref8}.

The ordinary Euler polynomials are defined by 
\begin{equation}\label{eq4}
\frac{2}{e^t+1}e^{xt}=\sum_{n=0}^\infty E_n^*(x)\frac{t^n}{n!},\quad
(\text{see \cite{ref1,ref7,ref10,ref11}}).
\end{equation}
When $x=0$, $E_n^*=E_n^*(0)$ are called the Euler numbers. 

Now, we define the type 2 Euler polynomials by 
\begin{equation}\label{eq5}
\frac{2}{e^{\frac{t}{2}}+e^{-\frac{t}{2}}}e^{xt} =\sum_{n=0}^\infty E_n(x)\frac{t^n}{n!},\quad (\text{see \cite{ref7,ref8,ref9}}).
\end{equation}
For $x=0$, $E_n=E_n(0)$ are called the type 2 Euler numbers so that they are given by
\begin{equation}\label{eq6}
\frac{2}{e^{\frac{t}{2}}+e^{-\frac{t}{2}}}=\sech\frac{t}{2}=\sum_{n=0}^\infty
E_n\frac{t^n}{n!}.
\end{equation}
Again, the type 2 Euler polynomials and numbers are slightly differently defined in \cite{ref8}.
From \eqref{eq4} and \eqref{eq6}, we note that
\begin{equation*}
E_n^* \left(\frac{1}{2}\right)=E_n,\quad (n\geq0),\quad (\text{see \cite{ref8}}).
\end{equation*}

Let $f$ be a uniformly differentiable function on  $\mathbb{Z}_p$. The bosonic (also called Volkenborn)  $p$-adic integral on $\mathbb{Z}_p$ is defined by
\begin{equation}\label{eq7}
I_1(f)=\intz f(x)
d\mu_1(x)=\lim_{N\to\infty}\frac{1}{p^N}\sum_{x=0}^{p^N-1}f(x),\quad
(\text{see \cite{ref10,ref11,ref13}}).
\end{equation}
From \eqref{eq7}, we note that
\begin{equation}\label{eq8}
I_1(f_1)-I_1(f)=f'(0),
\end{equation}
where $f_1(x)=f(x+1)$, and $f'(0)=\frac{df}{dx}\Big|_{x=0}$.

The fermionic $p$-adic integral on  $\mathbb{Z}_p$ was introduced by Kim as
\begin{equation}\label{eq9}
I_{-1}(f)=\intz f(x) d\mu_{-1}(x)=\lim_{N\to\infty}
\sum_{x=0}^{p^N-1}(-1)^xf(x),\quad (\text{see \cite{ref10,ref11}}).
\end{equation}
By \eqref{eq9}, we easily get
\begin{equation}\label{eq10}
I_{-1}(f_1)+I_{-1}(f)=2f(0).
\end{equation}

For $\lambda\in \mathbb{R}$, the degenerate exponential function is defined by
\begin{equation}\label{eq11}
e^x_{\lambda}(t)=(1+\lambda t)^{\frac{x}{\lambda}},\quad (\text{see \cite{ref3,ref4,ref7,ref8,ref9,ref12}}).
\end{equation}
Note that $\lim_{\lambda\to 0}e^x_{\lambda}(t)=e^{xt}$.
From \eqref{eq11} we have
\begin{equation}\label{eq12}
e^x_{\lambda}(t)=(1+\lambda t)^{\frac{x}{\lambda}}=\sum_{k=0}^\infty (x)_{k,\lambda}\frac{t^k}{k!},
\end{equation}
where $(x)_{k,\lambda}=x(x-\lambda)(x-2\lambda)\cdots (x-(k-1)\lambda),~(k \geq 1)$, and $ (x)_{0,\lambda}=1$.

\medskip
\section{Some identities of special polynomials arising from $p$-adic integrals
on $\mathbb{Z}_p$}\label{sec2}

From \eqref{eq8}, we note that
\begin{equation}\label{eq13}
\begin{aligned}
\int_{\mathbb{Z}_p}e^{(x+y+\frac{1}{2})t}d\mu_1(y)&=\frac{t}{e^{\frac{t}{2}}-e^{-\frac{t}{2}}}e^{xt}\\
&=\frac{t}{2}\csch\frac{t}{2}e^{xt}\\
&=\sum_{n=0}^\infty b_n(x) \frac{t^n}{n!}.
 \end{aligned}
\end{equation}
On the other hand, we have
\begin{equation}\label{eq14}
\begin{aligned}
\int_{\mathbb{Z}_p}e^{(x+y+\frac{1}{2})t}d\mu_1(x)&= \sum_{n=0}^\infty \int_{\mathbb{Z}_p}(x+y+\frac{1}{2})^n d\mu_1 (x)\frac{t^n}{n!}.
 \end{aligned}
\end{equation}
Therefore, by \eqref{eq13} and \eqref{eq14}, we obtain the following lemma.
\begin{lemma}
For $n \ge 0$, we have
\begin{equation*}
\intz \left( x+y+\frac{1}{2} \right)^n d \mu_1 (x) = b_n(x).
\end{equation*}
\end{lemma}

\medskip

By \eqref{eq7}, we get
\begin{equation} \label{eq15}
\begin{split}
\intz f(x) d \mu_1 (x) & = \lim_{N \rightarrow \infty} \frac{1}{p^N}
\sum_{x=0}^{p^N -1} f(x) = \lim_{N \rightarrow \infty}
\frac{1}{dp^N} \sum_{x=0}^{dp^N -1} f(x) \\
&= \frac{1}{d} \sum_{a=0}^{d -1} \lim_{N \rightarrow \infty}
\frac{1}{p^N} \sum_{x=0}^{p^N -1} f(a+xd) = \frac{1}{d}
\sum_{a=0}^{d -1} \intz f(a+ xd) d \mu_1 (x),
\end{split}
\end{equation}
where $d$ is a positive integer.

Therefore, by \eqref{eq15}, we obtain the following lemma.

\begin{lemma}
For $d  \in \mathbb{N}$, we have
\begin{equation*}
\intz f(x) d \mu_1 (x) = \frac{1}{d} \sum_{a=0}^{d -1} \intz f(a+
xd) d \mu_1 (x).
\end{equation*}
\end{lemma}

\medskip

Applying Lemma 2.2 to $f(x)= e^{(x+y+1/2)t}$, we have
\begin{equation}\label{eq16}
\begin{split}
\intz e^{(x+ y+ 1/2)t}  \dm (y) & =  \frac{1}{d} \sum_{a=0}^{d -1}
\intz e^{(x+a+dy+ 1/2)t} \dm (y)\\
&=  \frac{1}{d} \sum_{a=0}^{d -1} \intz e^{d(y+ \frac{1}{d}(x+a+
\frac{1-d}{2})+ 1/2)t} \dm (y).
\end{split}
\end{equation}

Thus, by \eqref{eq16}, we get

\begin{equation}\label{eq17}
\begin{split}
\sum_{n=0}^{\infty} & \intz \left(x+y+ \frac{1}{2} \right)^n \dm
(y) \frac{t^n}{n!} \\
&= \sum_{n=0}^{\infty}  d^{n-1} \sum_{a=0}^{d-1} \intz \left( y+
\frac{1}{d}\left( a+ x + \frac{1-d}{2} \right)+ 1/2 \right)^n \dm
(y) \frac{t^n}{n!}.
\end{split}
\end{equation}

By comparing the coefficients on both sides of \eqref{eq17}, we
get

\begin{equation}\label{eq18}
\intz \left(x+y+\frac{1}{2} \right)^n \dm (y)  =
 d^{n-1} \sum_{a=0}^{d-1} \intz \left( y+
\frac{1}{d} \left( a+ x + \frac{1-d}{2} \right)+ 1/2 \right)^n \dm
(y)
\end{equation}

By Lemma 2.1 and \eqref{eq18}, we get

\begin{equation}\label{eq19}
b_n (x) = d^{n-1} \sum_{a=0}^{d-1} b_n \left( \frac{x+a + \frac{1}{2}(1-d)}{d}
 \right),~ ( n \ge 0),
\end{equation}
where $d$ is a positive integer.

\begin{theorem} For $d \in \mathbb{N}$ and $n \in \mathbb{N} \cup
\{ 0 \}$, we have
\begin{equation*}
b_n (x) = d^{n-1} \sum_{a=0}^{d-1} b_n \left( \frac{x+a + \frac{1}{2}(1-d)}{d}
 \right).
\end{equation*}
\end{theorem}

\medskip

For $r \in \mathbb{N}$, we consider the multivariate $p$-adic
integral on $\mathbb{Z}_p$ as follows:

\begin{equation} \label{eq20}
\begin{split}
\intz & \cdots \intz e^{(x_1 + x_2 + \cdots + x_r + r/2)t} \dm (x_1)
\dm (x_2) \cdots \dm (x_r ) \\
&= \left( \frac{t}{ e^{t/2} - e^{-t/2}} \right)^r =
\left( \frac{t}{2} \csch \frac{t}{2} \right)^r.
\end{split}
\end{equation}

Now, we define the type 2 Bernoulli numbers of order $r$ by

\begin{equation} \label{eq21}
\left( \frac{t}{ e^{t/2} - e^{-t/2}} \right)^r =
\left( \frac{t}{2} \csch \frac{t}{2} \right)^r =
\sum_{n=0}^{\infty} b_n^{(r)} \frac{t^n}{n!}.
\end{equation}

By \eqref{eq20} and \eqref{eq21}, we see that
\begin{equation} \label{eq22}
\intz \cdots \intz \left( x_1 + x_2 +  \cdots + x_r + \frac{r}{2}
 \right)^n \dm (x_1 ) \cdots \dm (x_r ) = b_n^{(r)},~ ( n \ge 0).
\end{equation}

On the other hand,
\begin{equation} \label{eq23}
\begin{split}
\intz & \cdots \intz \left( x_1 + x_2 +  \cdots + x_r + \frac{r}{2}
 \right)^n \dm (x_1 ) \cdots \dm (x_r ) \\
  &= \sum_{i_1 + i_2 + \cdots + i_r =n \atop i_1, i_2 , \cdots , i_r
\ge 0 } \binom{n}{i_1 , i_2 , \cdots , i_r}
   \intz \left( x_1 + \frac{1}{2} \right)^{i_1} \dm  (x_1 ) \cdots
    \intz \left( x_r + \frac{1}{2} \right)^{i_r} \dm  (x_r ) \\
&= \sum_{i_1 + i_2 + \cdots + i_r =n \atop i_1, i_2 , \cdots , i_r
\ge 0 } \binom{n}{i_1 , i_2 , \cdots , i_r}
   b_{i_1} b_{i_2} \cdots b_{i_r}. \\
\end{split}
\end{equation}

Therefore, by \eqref{eq22} and \eqref{eq23}, we obtain the following
theorem.

\begin{theorem} For $n \ge 0, r \in \mathbb{N}$, we have
\begin{equation*}
b_n^{(r)} = \sum_{i_1 + i_2 + \cdots + i_r =n \atop i_1, i_2 ,
\cdots , i_r \ge 0 } \binom{n}{i_1 , i_2 , \cdots , i_r}
   b_{i_1} b_{i_2} \cdots b_{i_r}.
\end{equation*}
\end{theorem}

\medskip

From \eqref{eq21}, we have
\begin{equation}\label{eq24}
\begin{split}
t^r & = \sum_{l=0}^{\infty} b_l^{(r)} \frac{t^l}{l!} 
        \left( e^{\frac{t}{2}} - e^{-\frac{t}{2}} \right)^r =
        \sum_{l=0}^{\infty} b_l^{(r)} \frac{t^l}{l!}
        r! \sum_{m=r}^{\infty} T(m,r) \frac{t^m}{m!}\\
     &= \sum_{n=r}^{\infty} r! \sum_{m=r}^n \binom{n}{m} T(m, r) b_{n-m}^{(r)} \frac{t^n}{n!},
\end{split}
\end{equation}
where $T(m,r)$ are the central factorial numbers of the second kind.

Therefore, by \eqref{eq24}, we obtain the following theorem.

\begin{theorem}
For $n, r \in \mathbb{N} \cup \{ 0 \} $ with $ n \ge r $, we have
\begin{equation*}
\sum_{m=r}^n \binom{n}{m} T(m, r) b_{n-m}^{(r)} =
     \left\{ \begin{matrix} 1, & \text{if} & n=r, \\
                         0,  & \text{if} & n > r,
           \end{matrix} \right.
\end{equation*}
where $T(m,r)$ are the central factorial number of the second kind.
\end{theorem}

\medskip

From Lemma 2.1, we note that

\begin{equation}\label{eq25}
\begin{split}
b_n{(x)} &= \intz \left( y+ x + \frac{1}{2} \right)^n \dm (y)  =
\sum_{l=0}^n \binom{n}{l} x^{n-l} \intz \left( y+ \frac{1}{2}
\right)^l \dm (y)  \\
  &= \sum_{l=0}^{n} \binom{n}{l} x^{n-l} b_l .
\end{split}
\end{equation}

By \eqref{eq25}, we get
\begin{equation} \label{eq26}
b_n (x) = \sum_{l=0}^n \binom{n}{l} x^{n-l} b_l.
\end{equation}

Now, we observe that
\begin{equation} \label{eq27}
\begin{split}
\sum_{k=0}^{n-1} e^{\left(k+ \frac{1}{2} \right)t} &=
    e^{\frac{1}{2} t} \sum_{k=0}^{n-1} e^{kt} = \frac{1}{
    e^{\frac{t}{2}} - e^{-\frac{t}{2}}} (e^{nt}-1) \\
    &= \left( \frac{t}{
    e^{\frac{t}{2}} - e^{-\frac{t}{2}}} e^{nt} - \frac{t}{
    e^{\frac{t}{2}} - e^{-\frac{t}{2}}}  \right) \frac{1}{t} \\
    &= \frac{1}{t}  \sum_{m=0}^{\infty} (b_m (n) - b_m )
    \frac{t^m}{m!}  = \sum_{m=0}^{\infty} \frac{(b_{m+1} (n) -
    b_{m+1} )}{m+1} \frac{t^m}{m!}.
\end{split}
\end{equation}
On the other hand,
\begin{equation} \label{eq28}
\sum_{k=0}^{n-1} e^{ \left( k+ \frac{1}{2} \right) t} =
\sum_{m=0}^{\infty} \sum_{k=0}^{n-1} \left( k+ \frac{1}{2}
\right)^m  \frac{t^m}{m!}.
\end{equation}

By \eqref{eq27} and \eqref{eq28}, we get
\begin{equation} \label{eq29}
 \sum_{k=0}^{n-1} \left( 2k+ 1
\right)^m = 2^m \left( \frac{b_{m+1}(n) - b_{m+1} }{m+1} \right).
\end{equation}

Therefore, by \eqref{eq29}, and interchanging $m$ and $n$, we obtain the following theorem.

\begin{theorem}
For $m  \in \mathbb{N}$ and $n \in \mathbb{N} \cup \{ 0 \} $, we
have
\begin{equation*}
1^n + 3^n + \cdots + (2m-1)^n = 2^n \left( \frac{b_{n+1}(m) -
b_{n+1} }{n+1} \right).
\end{equation*}
\end{theorem}

\medskip

We define the fully degenerate type 2 Bernoulli polynomials by 

\begin{equation} \label{eq30}
\frac{1}{\lambda} \left( \frac{ \log(1+ \lambda
t)}{e_{\lambda}^{1/2} (t)- e_{\lambda}^{-1/2}(t)} \right)
e_{\lambda}^x (t) = \sum_{n=0}^{\infty} B_{n,\lambda}(x)
\frac{t^n}{n!}.
\end{equation}
When $x=0, B_{n, \lambda} =  B_{n, \lambda}(0)$ are called the fully
degenerate type 2 Bernoulli numbers.

We note that
\begin{equation} \label{eq31}
\begin{split}
\intz e_{\lambda}^{ x+ y + 1/2} (t) \dm (y) & = \frac{ \log ( 1+
\lambda t)}{ \lambda} \cdot \frac{1}{e_{\lambda}^{1/2} (t)-
e_{\lambda}^{-1/2}(t)} e_{\lambda}^x (t)  \\
&= \sum_{n=0}^{ \infty} B_{n, \lambda} (x) \frac{t^n}{n!}.
\end{split}
\end{equation}
Thus, by \eqref{eq31} and \eqref{eq12} we obtain
\begin{equation}\label{eq32}
\intz \left(x+y+ \frac{1}{2} \right)_{n, \lambda} \dm (y) = B_{n,
 \lambda}(x).
\end{equation}

As is known, the degenerate Stirling numbers of the first kind are
defined by
\begin{equation} \label{eq33}
(x)_{n, \lambda} = \sum_{l=0}^n S_{1,\lambda} (n, l) x^l,~ ( n \ge
0).
\end{equation}
By \eqref{eq32}, \eqref{eq33} and Lemma 2.1, we have
\begin{equation}\label{eq34}
B_{n,\lambda}(x)=\sum_{l=0}^n S_{1,\lambda} (n, l)b_l(x).
\end{equation}
Also, from \eqref{eq12} and \eqref{eq31} we observe that
\begin{equation} \label{eq35}
\begin{split}
\intz e_{\lambda}^{ x+ y + 1/2} (t) \dm (y) & = e_{\lambda}^x (t)
\intz e_{\lambda}^{y + 1/2} (t) \dm (y) \\
   &= \sum_{l=0}^{\infty} (x)_{l, \lambda} \frac{t^l}{l!}
   \sum_{m=0}^{\infty} B_{m,\lambda} \frac{t^m}{m!} \\
   &= \sum_{n=0}^{\infty}\sum_{m=0}^{n}
     \binom{n}{m}B_{m,\lambda} (x)_{n-m, \lambda}
   \frac{t^n}{n!}.
\end{split}
\end{equation}

Therefore, from \eqref{eq32}, \eqref{eq34} and \eqref{eq35},  we have the following theorem.

\begin{theorem} For $n\ge 0$, we have
\begin{equation*}
B_{n,\lambda}(x)=\intz \left(x+y+ \frac{1}{2} \right)_{n, \lambda} \dm (y) =
\sum_{l=0}^n S_{1,\lambda} (n, l)b_l(x)=\sum_{m=0}^{n}
\binom{n}{m}B_{m,\lambda} (x)_{n-m, \lambda}.
\end{equation*}
\end{theorem}

\medskip

As is known, the degenerate Carlitz type 2 Bernoulli polynomials are
defined by
\begin{equation} \label{eq36}
\frac{t}{e_{\lambda}^{\frac{1}{2}}(t) -
e_{\lambda}^{-\frac{1}{2}}(t)} e_{\lambda}^x (t) =
\sum_{n=0}^{\infty} b_{n, \lambda} (x) \frac{t^n}{n!}.
\end{equation}
When $x=0, b_{n, \lambda} = b_{n, \lambda}(0), (n \ge 0)$, are
called the degenerate Carlitz type 2 Bernoulli numbers.

It is well known that the Daehee numbers, denoted by $d_n$, are
defined by 
\begin{equation} \label{eq37}
\frac{ \log (1+t)}{t} = \sum_{n=0}^{\infty} d_n \frac{t^n}{n!},\quad (\text{see \cite{ref6, ref16}}).
\end{equation}

Now, from \eqref{eq31},\eqref{eq36} and \eqref{eq37}, we observe that
\begin{equation} \label{eq38}
\begin{split}
\sum_{n=0}^{\infty}B_{n,\lambda}\frac{t^n}{n!}=\intz e_{\lambda}^{ x+ 1/2}(t) \dm (x) &= \frac{ \log(1+ \lambda
t)}{ \lambda t} \frac{t}{  e_{\lambda}^{1/2}(t) -
e_{\lambda}^{-1/2}(t)} \\
&= \sum_{l=0}^{\infty} \lambda^l d_l \frac{t^l}{l!}
\sum_{m=0}^{\infty} b_{m,\lambda} \frac{t^m}{m!}\\
&= \sum_{n=0}^{\infty}\sum_{l=0}^{n} \binom{n}{l} \lambda^l
d_l  b_{n-l,\lambda}\frac{t^n}{n!}.
\end{split}
\end{equation}

Therefore, by \eqref{eq38} and \eqref{eq12}, we obtain the following theorem.

\begin{theorem} For $n \ge 0$, we have
\begin{equation*}
B_{n,\lambda} =\intz \left( x + \frac{1}{2} \right)_{n, \lambda} \dm (x)=
\sum_{l=0}^{n} \binom{n}{l} \lambda^l d_l b_{n-l, \lambda}.
\end{equation*}
\end{theorem}

\medskip

For $n \in \mathbb{N}$, by \eqref{eq8}, we easily get

\begin{equation} \label{eq39}
\intz f(x+m) \dm (x) =  \sum_{l=0}^{m-1} f ' (x)  + \intz f(x) \dm
(x).
\end{equation}

By applying \eqref{eq39} to $f(x)=e_{\lambda}^{x+\frac{1}{2}}(t)$, we get
\begin{equation} \label{eq40}
\frac{1}{ e_{\lambda}^{1/2}(t)  - e_{\lambda}^{-1/2}(t) }
e_{\lambda}^m (t) - \frac{1}{ e_{\lambda}^{1/2}(t)  -
e_{\lambda}^{-1/2}(t) } = e_{\lambda}^{1/2}(t) \sum_{l=0}^{m-1}
e_{\lambda}^{l}(t).
\end{equation}

From \eqref{eq40}, we derive the following equation.
\begin{equation} \label{eq41}
\frac{1}{t} \sum_{n=0}^{\infty} ( b_{n, \lambda} (m)  - b_{n,
\lambda}) \frac{t^n}{n!} = \sum_{n=0}^{\infty} \left(
\sum_{l=0}^{m-1} \left(l+ \frac{1}{2} \right)_{n, \lambda} \right)
\frac{t^n}{n!}.
\end{equation}

By \eqref{eq41}, we get
\begin{equation} \label{eq42}
\sum_{n=0}^{\infty} \left( \frac{b_{n+1,\lambda}{(m)} -b_{n+1,
\lambda}}{ n+1} \right) \frac{t^n}{n!} = \sum_{n=0}^{\infty} \left(
\frac{1}{2^n} \sum_{l=0}^{m-1} (2l+1)_{n, 2\lambda}\right)
\frac{t^n}{n!}.
\end{equation}

Therefore, by \eqref{eq42}, we obtain the following theorem.

\begin{theorem}
For $n \ge 0,  m \in \mathbb{N}$, we have
\begin{equation*}
\frac{2^n}{ n+1} ( b_{n+1,\lambda}{(m)} -b_{n+1, \lambda}) =
\sum_{l=0}^{m-1} (2l+1)_{n, 2\lambda}.
\end{equation*}
\end{theorem}

\medskip

From \eqref{eq10}, we observe that
\begin{equation} \label{eq43}
\intz e^{t \left( x+ y+\frac{1}{2} \right)} \dmm (y) = \frac{2}{
e^{t/2} + e^{-t/2}}e^{xt} = \sech\frac{t}{2}e^{xt}=\sum_{n=0}^{\infty} E_n (x)
\frac{t^n}{n!}.
\end{equation}

Thus from \eqref{eq43} and \eqref{eq12}, we have the following lemma.

\begin{lemma}
For $n \ge 0$, we have
\begin{equation*}
\intz \left( x+y+\frac{1}{2} \right)^n \dmm (y) = E_n(x).
\end{equation*}
\end{lemma}

From Lemma 2.10, we have
\begin{equation}\label{eq44}
\begin{split}
E_n (x) = \intz \left( x+y+ \frac{1}{2} \right)^n \dmm(y) &=
\sum_{l=0}^{n} \binom{n}{l} x^{n-l} \intz \left( y +
\frac{1}{2} \right)^l \dmm (y) \\
&= \sum_{l=0}^{n} \binom{n}{l}x^{n-l} E_l ,~ (n \ge 0 ).
\end{split}
\end{equation}

Let $d \in \mathbb{N}$ with $ d \equiv 1 (mod 2)$. Then, by
\eqref{eq10}, we get

\begin{equation} \label{eq45}
\intz f(x+d) \dmm (x) + \intz f(x) \dmm (x) = 2 \sum_{l=0}^{d-1}
(-1)^{l} f(l).
\end{equation}

Let us take $f(x) = e^{(x+1/2)t}$. Then, by \eqref{eq45}, we get
\begin{equation} \label{eq46}
e^{mt} \intz e^{(x+1/2)t} \dmm (x) + \intz e^{(x+ 1/2)t} \dmm (x) =
2 \sum_{l=0}^{m-1} (-1)^l e^{(l + 1/2)t}.
\end{equation}

From \eqref{eq46}, we have
\begin{equation}\label{eq47}
\frac{2}{e^{t/2} + e^{-t/2}} e^{mt} + \frac{2}{e^{t/2} + e^{-t/2}} =
2\sum_{l=0}^{m-1} (-1)^l e^{(l+1/2)t}.
\end{equation}

By \eqref{eq5} and \eqref{eq47}, we get
\begin{equation}\label{eq48}
\sum_{n=0}^{\infty} ( E_n (m) + E_n ) \frac{t^n}{n!} =
\sum_{n=0}^{\infty} \left( 2 \sum_{l=0}^{m-1} (-1)^l \left( l+
\frac{1}{2}\right)^n \right) \frac{t^n}{n!}.
\end{equation}

Therefore, by \eqref{eq48}, we obtain the following theorem.

\begin{theorem}
For $m \in \mathbb{N}$ with $m \equiv 1 (mod 2),~ n  \in
\mathbb{N} \cup \{ 0 \}$, we have
\begin{equation*}
2^{n-1} ( E_n (m) + E_n )  = \sum_{l=0}^{m-1} (-1)^l (2l +1)^n.
\end{equation*}
\end{theorem}

\medskip

The following lemma can be easily shown.

\begin{lemma}
\begin{equation*}
\intz f(x) \dmm (x) = \sum_{a=0}^{d-1} (-1)^a \intz f( a+ dx ) \dmm
(x),
\end{equation*}
\end{lemma}
where $ d \in \mathbb{N} $ with $ d \equiv 1 ( \mod 2)$.

\medskip

Let us apply Lemma 2.12 to $f(y) = (x+ y+ 1/2)^n$. Then we have

\begin{equation} \label{eq49}
\begin{split}
\intz \left( x+ y + \frac{1}{2} \right)^n \dmm(y) & =
\sum_{a=0}^{d-1} (-1)^a \intz \left( x + a + dy + \frac{1}{2}
\right)^n  \dmm  (y) \\
 &= d^n \sum_{a=0}^{d-1} (-1)^a \intz \left( \frac{x+ a+ \frac{1}{2}(1-d)}{d} + y + \frac{1}{2}
   \right)^n  \dmm (y).
\end{split}
\end{equation}

Therefore, by \eqref{eq49}, we have the following theorem.

\begin{theorem}
For $d \in \mathbb{N}$ with $d \equiv 1 (mod 2 ),~ n \in \mathbb{N}
\cup \{ 0 \} $, we have
\begin{equation*}
E_n (x) = d^n \sum_{a=0}^{d-1} (-1)^a E_n \left( \frac{x+a +
\frac{1}{2}(1-d)}{d} \right).
\end{equation*}
\end{theorem}

\medskip

For $r \in \mathbb{N}$, let us consider the following fermionic
$p$-adic integral on $\mathbb{Z}_p$.
\begin{equation} \label{eq50}
\begin{split}
\intz \intz & \cdots \intz e^{ \left(x_1 + x_2 + \cdots + x_r +
\frac{r}{2} \right)t} \dmm (x_1) \dmm (x_2) \cdots \dmm(x_r) \\
& = \left( \frac{2}{e^{\frac{t}{2}} + e^{-\frac{t}{2}}
}\right)^r = \left(\sech \frac{t}{2}\right)^r
\end{split}
\end{equation}

Let us define the type 2 Euler numbers of order $r$ by

\begin{equation} \label{eq51}
\left( \frac{2}{e^{\frac{t}{2}} + e^{-\frac{t}{2}}
}\right)^r = \left(\sech \frac{t}{2}\right)^r=
\sum_{n=0}^{\infty} E_n^{(r)} \frac{t^n}{n!}.
\end{equation}

From \eqref{eq50} and \eqref{eq51}, we have
\begin{equation} \label{eq52}
\intz \intz  \cdots \intz \left( x_1 + x_2 + \cdots + x_r +
\frac{r}{2} \right)^n \dmm (x_1) \dmm (x_2) \cdots \dmm(x_r) =
E_n^{(r)},~(n \ge 0).
\end{equation}

On the other hand,
\begin{equation} \label{eq53}
\begin{split}
\intz \intz & \cdots \intz \left( x_1 + x_2 + \cdots + x_r +
\frac{r}{2} \right)^n \dmm (x_1) \dmm (x_2) \cdots \dmm(x_r)\\
& = \sum_{ i_1 + i_2 + \cdots + i_r = n \atop i_1, i_2 , \cdots ,
i_r  \ge 0} \binom{n}{i_1 , \cdots , i_r } \intz \left( x_1 +
\frac{1}{2} \right)^{i_1} \dmm (x_1) \cdots \intz \left( x_r +
\frac{1}{2} \right)^{i_r} \dm (x_r) \\
&= \sum_{ i_1 + i_2 + \cdots + i_r = n \atop i_1, i_2 , \cdots , i_r
\ge 0 } \binom{n}{i_1 , \cdots , i_r } E_{i_1} E_{i_2} \cdots
E_{i_r}.
\end{split}
\end{equation}

Therefore, by \eqref{eq52} and \eqref{eq53}, we obtain the following theorem.

\begin{theorem}
For $n \ge 0$, we have
\begin{equation*}
E_n^{(r)} = \sum_{  i_1 + i_2 + \cdots + i_r = n \atop i_1, i_2 ,
\cdots , i_r  \ge 0 } \binom{n}{i_1 , \cdots , i_r } E_{i_1} E_{i_2}
\cdots E_{i_r}.
\end{equation*}
\end{theorem}

\medskip

From \eqref{eq51}, we have

\begin{equation} \label{eq54}
\begin{split}
2^r &= \sum_{l=0}^{\infty} E_l^{(r)}\frac{t^l}{l!}(e^{\frac{t}{2}} +e^{-\frac{t}{2}} )^r  \\
&= \sum_{l=0}^{\infty} E_l^{(r)}\frac{t^l}{l!}
\sum_{j=0}^{r} \binom{r}{j}e^{\left(j-\frac{r}{2} \right) t} \\
&= \sum_{l=0}^{\infty} E_l^{(r)}\frac{t^l}{l!}
\sum_{m=0}^{\infty}  \sum_{j=0}^{r} \binom{r}{j}\left(j-\frac{r}{2}
\right)^m \frac{t^m}{m!} \\
&= \sum_{n=0}^{\infty}\sum_{m=0}^n \sum_{j=0}^{r}
\binom{r}{j} \binom{n}{m} \left( j- \frac{r}{2} \right)^m
E_{n-m}^{(r)} \frac{t^n}{n!}.
\end{split}
\end{equation}

Comparing the coefficients on both sides of \eqref{eq54}, we
obtain the following theorem.

\begin{theorem}
For $n \ge 0$, we have
\begin{equation*}
\sum_{m=0}^n \sum_{j=0}^{r} \binom{r}{j} \binom{n}{m} \left( j-
\frac{r}{2} \right)^m E_{n-m}^{(r)}  = \left\{ \begin{matrix} 2^r, &
\text{if}& n=0, \\  0, & \text{if}& n > 0 .
\end{matrix} \right.
\end{equation*}
\end{theorem}

\medskip

We define the degenerate type 2 Euler polynomials by
\begin{equation} \label{eq55}
\frac{2}{e_{\lambda}^{1/2}(t) + e_{\lambda}^{-1/2}(t)}
e_{\lambda}^{x}(t)
= \sum_{n=0}^{\infty} E_{n, \lambda}(x) \frac{t^n}{n!}.
\end{equation}
When $x=0$, $E_{n,\lambda}=E_{n,\lambda}(0)$ are called the degenerate type 2 Euler numbers.

From \eqref{eq10}, we can derive the following equation.
\begin{equation} \label{eq56}
\begin{split}
\intz e_{\lambda}^{x+ y+ \frac{1}{2}} (t) \dmm (y) & =
\frac{2}{e_{\lambda}^{1/2}(t) + e_{\lambda}^{-1/2}(t)}
e_{\lambda}^{x}(t) \\
&= \sum_{n=0}^{\infty} E_{n, \lambda}(x) \frac{t^n}{n!}.
\end{split}
\end{equation}
By \eqref{eq56} and \eqref{eq12}, we get
\begin{equation} \label{eq57}
E_{n, \lambda}(x) = \intz \left( x+ y + \frac{1}{2} \right)_{n,
\lambda} \dmm (y),~ ( n \ge 0 ).
\end{equation}

By \eqref{eq57}, \eqref{eq33} and Lemma 2.10, we get
\begin{equation} \label{eq58}
E_{n,\lambda} (x) = \sum_{l=0}^n S_{1,\lambda} (n, l) E_l (x).
\end{equation}

Also, from \eqref{eq12} and \eqref{eq56}, we observe that
\begin{equation} \label{eq59}
\begin{split}
\intz e_{\lambda}^{ x+ y + 1/2} (t) \dmm (y) & = e_{\lambda}^x (t)
\intz e_{\lambda}^{y + 1/2} (t) \dmm (y) \\
&= \sum_{l=0}^{\infty} (x)_{l, \lambda} \frac{t^l}{l!}
\sum_{m=0}^{\infty} E_{m,\lambda} \frac{t^m}{m!} \\
&= \sum_{n=0}^{\infty}\sum_{m=0}^{n}
\binom{n}{m}E_{m,\lambda} (x)_{n-m, \lambda}
\frac{t^n}{n!}.
\end{split}
\end{equation}

Therefore, by \eqref{eq57}-\eqref{eq59}, we obtain the following theorem.

\begin{theorem}
For $n \ge 0$, we have
\begin{equation*}
E_{n,\lambda}(x)= \intz \left(x+ y+ \frac{1}{2}\right)_{n,\lambda} \dmm (y) = \sum_{l=0}^n  S_{1,\lambda} (n, l) E_l (x)=\sum_{m=0}^{n}\binom{n}{m}E_{m,\lambda} (x)_{n-m, \lambda}.
\end{equation*}
\end{theorem}

\medskip

For $m \in \mathbb{N}$ with $ m \equiv 1 (mod 2)$, from \eqref{eq45} we have
\begin{equation} \label{eq60}
\intz e_{\lambda}^{m+x+ 1/2} (t) \dmm (x) +  \intz
e_{\lambda}^{x+1/2} (t) \dmm (x) = 2 \sum_{l=0}^{m-1} (-1)^l
e_{\lambda}^{l+ 1/2}(t).
\end{equation}

From \eqref{eq60}, we have
\begin{equation} \label{eq61}
\begin{split}
\sum_{n=0}^{\infty} \left( E_{n, \lambda} (m) + E_{n, \lambda}
\right) \frac{t^n}{n!} & = 2 \sum_{n=0}^{\infty} \sum_{l=0}^{m-1} (-1)^l
\left(l + \frac{1}{2} \right)_{n, \lambda}
\frac{t^n}{n!} \\
&= \sum_{n=0}^{\infty} \left(  \frac{1}{2} \right)^{n-1}
\sum_{l=0}^{m-1} (-1)^l (2l+1)_{n, 2\lambda} \frac{t^n}{n!}.
\end{split}
\end{equation}

Therefore, by \eqref{eq61}, we obtain the following theorem.

\begin{theorem}
For $n \ge 0, m \in \mathbb{N}$ with $m \equiv 1 (mod 2 )$, we
have
\begin{equation*}
2^{n-1} \left( E_{n, \lambda} (m) + E_{n,  \lambda} \right) =
\sum_{l=0}^{m-1} (-1)^l (2l+1)_{n, 2 \lambda}.
\end{equation*}

\end{theorem}

\medskip

For $r \in \mathbb{N}$, we have
\begin{equation} \label{eq62}
\begin{split}
\intz  & \cdots \intz e_{\lambda}^{x_1 + \cdots+ x_r + r/2} (t) \dmm
(x_1) \dmm (x_2 ) \cdots \dmm (x_r) \\
&=  \left( \frac{2}{e_{\lambda}^{1/2}(t)+
e_{\lambda}^{-1/2}(t)}  \right)^r
\end{split}
\end{equation}

Now, we define the degenerate type 2 Euler numbers of order $r$ which are given
by
\begin{equation} \label{eq63}
\left( \frac{2}{e_{\lambda}^{1/2}(t)+
e_{\lambda}^{-1/2}(t)}  \right)^r = \sum_{n=0}^{\infty} E_{n, \lambda}^{(r)}
\frac{t^n}{n!}.
\end{equation}

By \eqref{eq62}, \eqref{eq63} and \eqref{eq12}, we get

\begin{equation*}
\intz \cdots \intz \left( x_1 + x_2 + \cdots + x_r + \frac{r}{2}
\right)_{n,\lambda} \dmm (x_1) \dmm(x_2 ) \cdots \dmm (x_r)=E_{n, \lambda}^{(r)},~ (n \ge 0).
\end{equation*}

\section{Conclusion}\label{sec3}

\indent In recent years, studying degenerate versions of various special polynomials and numbers have attracted many mathematicians and been carried out by several different methods like generating functions, combinatorial approaches, umbral calculus, $p$-adic analysis and differential equations.
In this paper, we introduced degenerate type 2 Bernoulli polynomials, fully degenerate type 2 Bernoulli polynomials and degenerate type 2 Euler polynomials, and their corresponding numbers, as degenerate and type 2 versions of Bernoulli and Euler numbers. We investigated those polynomials and numbers by means of  bosonic and fermionic $p$-adic integrals and derived some identities, distribution relations, Witt type formulas and analogues for the Bernoulli's interpretation of powers of the first $m$ positive integers in terms of Bernoulli polynomials. In more detail, our main results are as follows.\\
\indent As to the analogues for the Bernoulli's interpretation of power sums, in Theorem 2.6 we expressed powers of the first $m$ odd integers in terms of type 2 Bernoulli polynomials $b_n(x)$, in Theorem 2.11  alternating sum of powers of the first $m$ odd integers in terms of type 2 Euler polynomials $E_n(x)$, in Theorem 2.9 sum of the values of the generalized falling factorials at the first $m$ odd positive integers in terms of degenerate Carlitz type 2 Bernoulli polynomials $b_{n,\lambda}(x)$, and in Theorem 2.17 alternating sum of the values of the generalized falling factorials at the first $m$ odd positive integers in terms of degenerate type 2 Euler polynomials $E_{n,\lambda}(x)$. Witt type formulas were obtained for $b_n(x), B_{n,\lambda}(x), E_n(x)$, and $E_{n,\lambda}(x)$, respectively in Lemma 2.1, Theorem 2.7, Lemma 2.10 and Theorem 2.16. 
Distribution relations were derived for  $b_n(x)$, and $E_n(x)$, respectively in Theorem 2.3 and Theorem 2.13. \\
\indent As one of our future projects, we would like to continue to do researches on degenerate versions of various special numbers and polynomials, and find many applications of them in mathematics, sciences and engineering.


\begin{thebibliography}{10}

\bibitem{ref1}
S. Araci, M. Acikgoz,
\newblock A note on the Frobenius-Euler numbers and polynomials associated with Bernstein polynomials,
\newblock{\em Adv. Stud. Contmp. Math. (Kyungshang)} 22 (2012), no. 3, 399--406.

\bibitem{ref2}
A. Bayad, J. Chikhi,
\newblock Apostol-Euler polynomials and asymptotics for negative binomial reciproals,
\newblock {\em Adv. Stud. Contmp. Math. (Kyungshang)} 24 (2014), no. 1, 33--37.

\bibitem{ref3}
L. Carlitz,
\newblock Degenerate Stirling, Bernoulli and Eulerian numbers,
\newblock {\em Utilitas Math.} 15 (1979), 51--88.


\bibitem{ref4}
L. Carlitz,
\newblock A degenerate Staudt-Clausen theorem,
\newblock {\em Arch. Math, (Basel)} 7 (1956), 28--33.

\bibitem{ref5}
S. Gaboury, R. Tremblay, B.-J. Fug\`{e}re,
\newblock Some explicit formulas for cirtain new classes of Bernoulli, Euler and Genocchi polynomials,
\newblock {\em Proc. Jangjeon Math. Soc.} 17 (2014), no. 1, 115--123.

\bibitem{ref6}
G.-W. Jang, J. Kwon, J. G. Lee,
\newblock Some identities of degenerate Daehee numbers arising from nonlinear differential equation,
\newblock {\em Adv. Difference Equ.} (2017), 2017:206, 10 pp.

\bibitem{ref7}
G.-W. Jang, T. Kim,
\newblock A note on type 2 degenerate Euler and Bernoulli polynomials,
\newblock {\em Adv. Stud. Contmp. Math. (Kyungshang)} 29 (2019), no. 1, 147--159.

\bibitem{ref8}
T. Kim, D. S. Kim,
\newblock A note on type 2 Changhee and Daehee polynomials,
\newblock{\em Rev. R. Acad. Cienc. Exactas F\`{i}s, Nat. Ser. A Mat.
RACSAM(2019)} in press.

\bibitem{ref9}
T. Kim, D. S. Kim,
\newblock Degenerate central factorial numbers of the second kind,
\newblock {\em arXiv:1902.04360} [pdf, ps, other].


\bibitem{ref10}
T. Kim,
\newblock Some identities on the $q$-Euler polynomials of higher order and $q$-Stirling numbers
 by the fermionic $p$-adic integral on $\mathbb{Z}_p$,
\newblock{Russ. J. Math. Phys.} 16 (2009), no. 4, 484--491.


\bibitem{ref11}
T. Kim,
\newblock Symmetry of power sum polynomials and multivariate fermionic $p$-adic invariant integral
 on $\mathbb{Z}_p$,
\newblock{Russ. J. Math. Phys.}16 (2009), no. 1, 93--96.


\bibitem{ref12}
T. Kim, D. S. Kim,
\newblock Degenerate Laplace transform and degenerate gamma function,
\newblock{Russ. J. Math. Phys.} 24 (2017), no. 2, 241--248.

\bibitem{ref13}
L.-C. Jang, W.-J. Kim, Y. Simsek,
\newblock A study on the $p$-adic integral representation on $ \mathbb{Z}_p $ associated with
Bernstein and Bernoulli polynomials,
\newblock{Adv. Difference Equ.} 2010. Art. ID 163217, 6 pp.

\bibitem{ref14}
J. G. Lee, J. Kwon,
\newblock The modified degenerate $q$-Bernoulli polynomials arising from $p$-adic invariant integral on $\mathbb{Z}_p$,
\newblock{Adv. Difference Equ.} (2017), 2017:29, 9 pp.

\bibitem{ref15}
Q.-M. Luo,
\newblock Some recursion formulae and relations for Bernoulli numbers and Euler numbers of higher order,
\newblock{Adv. Stud. Contmp. Math. (Kyungshang)} 10 (2005), no. 1, 63--70.

\bibitem{ref16}
J.-W. Park, B. M. Kim, J. Kwon,
\newblock On a modified degenerate Daehee polynomials and numbers,
\newblock {\em J. Nonlinear Sci. Appl.} 10 (2017), no. 3, 1105--1115.

\bibitem{ref17}
Y. Simsek,
\newblock Generating functions of the twisted Bernoulli numbers and polynomials associated with teir interpolation functions,
\newblock {\em Adv. Stud. Contmp. Math. (Kyungshang)} 16 (2008), no. 2, 251--278.









\end{thebibliography}
\end{document}